\documentclass[12pt]{amsart}
\usepackage{times}

\begin{document}
\numberwithin{equation}{section}

\def\1#1{\overline{#1}}
\def\2#1{\widetilde{#1}}
\def\3#1{\widehat{#1}}
\def\4#1{\mathbb{#1}}
\def\5#1{\frak{#1}}
\def\6#1{{\mathcal{#1}}}

\newcommand{\de}{\partial}
\newcommand{\R}{\mathbb R}
\newcommand{\al}{\alpha}
\newcommand{\tr}{\widetilde{\rho}}
\newcommand{\tz}{\widetilde{\zeta}}
\newcommand{\tv}{\widetilde{\varphi}}
\newcommand{\hv}{\hat{\varphi}}
\newcommand{\tu}{\tilde{u}}
\newcommand{\tF}{\tilde{F}}
\newcommand{\debar}{\overline{\de}}
\newcommand{\Z}{\mathbb Z}
\newcommand{\C}{\mathbb C}
\newcommand{\Po}{\mathbb P}
\newcommand{\zbar}{\overline{z}}
\newcommand{\G}{\mathcal{G}}
\newcommand{\So}{\mathcal{S}}
\newcommand{\Ko}{\mathcal{K}}
\newcommand{\U}{\mathcal{U}}
\newcommand{\B}{\mathbb B}
\newcommand{\oB}{\overline{\mathbb B}}
\newcommand{\Cur}{\mathcal D}
\newcommand{\Dis}{\mathcal Dis}
\newcommand{\Levi}{\mathcal L}
\newcommand{\SP}{\mathcal SP}
\newcommand{\Sp}{\mathcal Q}
\newcommand{\Ma}{\mathcal M}
\newcommand{\Co}{\mathcal C}
\newcommand{\Eo}{\mathcal E}
\newcommand{\Hol}{{\sf Hol}(\mathbb H, \mathbb C)}
\newcommand{\Aut}{{\sf Aut}(\mathbb D)}
\newcommand{\D}{\mathbb D}
\newcommand{\oD}{\overline{\mathbb D}}
\newcommand{\oX}{\overline{X}}
\newcommand{\loc}{L^1_{\rm{loc}}}
\newcommand{\la}{\langle}
\newcommand{\ra}{\rangle}
\newcommand{\thh}{\tilde{h}}
\newcommand{\N}{\mathbb N}
\newcommand{\kd}{\kappa_D}
\newcommand{\Hr}{\mathbb H}
\newcommand{\ps}{{\sf Psh}}

\newcommand{\subh}{{\sf subh}}
\newcommand{\harm}{{\sf harm}}
\newcommand{\ph}{{\sf Ph}}
\newcommand{\tl}{\tilde{\lambda}}
\newcommand{\ts}{\tilde{\sigma}}
\newcommand{\tp}{\tilde{\phi}}
\newcommand{\defeq}{\mathrel{\mathop:}=}

\def\v{\varphi}
\def\Re{{\sf Re}\,}
\def\Im{{\sf Im}\,}
\def\Arg{{\sf Arg}\,}

\def\dist{{\rm dist}}
\def\const{{\rm const}}
\def\rk{{\rm rank\,}}
\def\id{{\sf id}}
\def\aut{{\sf aut}}
\def\Aut{{\sf Aut}}
\def\CR{{\rm CR}}
\def\GL{{\sf GL}}
\def\U{{\sf U}}

\def\la{\langle}
\def\ra{\rangle}

\newtheorem{theorem}{Theorem}[section]
\newtheorem*{maintheorem}{Main Theorem}
\newtheorem{lemma}[theorem]{Lemma}
\newtheorem{proposition}[theorem]{Proposition}
\newtheorem{corollary}[theorem]{Corollary}
\newtheorem{problem}[theorem]{Problem}

\theoremstyle{definition}
\newtheorem{definition}[theorem]{Definition}
\newtheorem{example}[theorem]{Example}

\theoremstyle{remark}
\newtheorem{remark}[theorem]{Remark}
\numberwithin{equation}{section}

\title{Valiron's construction in higher dimension}
\author[F. Bracci]{Filippo Bracci}
\address{F. Bracci: Dipartimento Di Matematica\\
Universit\`{a} di Roma \textquotedblleft Tor Vergata\textquotedblright\ \\
Via Della Ricerca Scientifica 1, 00133 \\
Roma, Italy} \email{fbracci@mat.uniroma2.it}
\thanks{}
\author[G. Gentili]{Graziano Gentili}\address{G.
Gentili: Dipartimento di Matematica "Ulisse Dini", Universit\`a
degli Studi di Firenze, Viale Morgagni 67/A, 50134 Firenze,
Italy}\email{gentili@math.unifi.it}
\author[P. Poggi-Corradini]{Pietro Poggi-Corradini}\address{Department of Mathematics, Cardwell Hall, Kansas State University,
Manhattan, KS 66506, USA.} \email{pietro@math.ksu.edu}

\thanks{Poggi-Corradini thanks the Department of Mathematics
  ``Ulisse Dini'' of the University of Florence for the hospitality,
  and the GNSAGA group for its financial support during the first stages of this project. }




\begin{abstract} We consider holomorphic self-maps $\v$ of the unit ball
  $\B^N$ in $\C^N$ ($N=1,2,3,...$). In the one-dimensional case, when
  $\v$ has no fixed points in $\D\defeq \B^1$ and is of hyperbolic
  type, there is a classical renormalization procedure due to Valiron
  which allows to semi-linearize the map $\phi$, and therefore, in
  this case, the
  dynamical properties of $\phi$ are well understood. In what follows,
we generalize the classical Valiron construction to higher
dimensions under some  weak assumptions on $\v$ at its
Denjoy-Wolff point. As a result, we construct a
semi-conjugation $\sigma$, which maps the ball into the right
half plane of $\C$, and solves the functional equation
$\sigma\circ \v=\lambda \sigma$, where $\lambda>1$ is the
(inverse of the) boundary dilation coefficient at the
Denjoy-Wolff point of $\v$.
\end{abstract}

\maketitle

\section{Introduction}

\subsection{The one-dimensional case}

Let $\v$ be a
holomorphic map on $\D$ with $\v(\D)\subset\D$. If $\v$ has no fixed points in $\D$, then by
the classical Wolff lemma (see, {\sl e.g.}, \cite{Abate}) there
exists a unique point $\tau\in\de\D$, called {\sl the Denjoy-Wolff
point of $\v$}, such that the sequence of iterates $\{\v^{\circ n}\}$ of
$\v$ converges uniformly on compacta to the constant map
$\zeta\mapsto \tau$, $\forall \zeta\in\D$. Also, by the classical
Julia-Wolff-Caratheodory theorem,  $\tau$ is a fixed point (as
nontangential limit) for $\v$ and the first derivative $\v'$
has nontangential limit $c\in (0,1]$ at $\tau$, moreover,
\[
c=\liminf_{\zeta\to\tau}\frac{1-|\v(\zeta)|}{1-|\zeta|}.
\]
The number $c$ is called the {\sl multiplier} of $\v$ or the
{\sl boundary dilatation coefficient} at $\tau$. The map $\v$
is called {\sl hyperbolic} if $c<1$ and {\sl parabolic} if
$c=1$.

Geometrically, one defines the horodisks $H(t)\defeq\{z\in\D:
|\tau-z|^2/(1-|z|^2)<1/t\}$,  which are disks in $\D$ internally tangent
to $\de \D$ at $\tau$, and which get smaller as $t$ gets larger.
Then the following mapping property holds: $\phi(H(t))\subset
H(t/c)$. In formulas:
\[
\frac{|\tau-\v(z)|^2}{1-|\v(z)|^2}\leq c\frac{|\tau-z|^2}{1-|z|^2},
\]
for every $z\in \D$.

 In 1931 G. Valiron
\cite{va1} (see also \cite{va2} and \cite{Br-Po}) proved that
if $\v$ is hyperbolic then there exists a nonconstant
holomorphic map $\theta:\D\to \Hr:=\{w\in \C: \Re w>0\}$ which
solves the so-called {\sl Schr\"oder equation}:
\begin{equation}\label{schroder}
\theta \circ \v= \frac{1}{c} \theta.
\end{equation}
Valiron constructs the map $\theta$ as follows.
First, in order to simplify notations, one can move to the
right half-plane $\Hr$ via the Cayley map
$C(\zeta)=(\tau+\zeta)/(\tau-\zeta)$, which takes
$\tau$ to $\infty$ and conjugates $\v$ to a self-map $\phi:=C\circ
\v \circ C^{-1}$ of $\Hr$, with Denjoy-Wolff point $\infty$
and multiplier $1/c$. Then, one considers the orbit
$x_n+iy_n:=\phi^{\circ n}(1)$ of the point $w=1$, and studies
the sequence of renormalized iterates:
\begin{equation}\label{valiron-method}
\sigma_n(w):= \frac{\phi^{\circ n}(w)}{x_n}.
\end{equation}
Valiron showed that $\{\sigma_n\}$ converges to a holomorphic
map $\sigma:\Hr\to \Hr$ such that $\sigma\circ
\phi=\frac{1}{c}\sigma$. Thus $\theta:=\sigma \circ C$ solves
\eqref{schroder}.

After Valiron's construction, Ch. Pommerenke \cite{pom},
\cite{pom2}, C. Cowen \cite{Co} and P. Bourdon and J. Shapiro
\cite{bs} exploited other constructions to solve
\eqref{schroder} (and the corresponding Abel's equation for the
parabolic case). In particular, Pommerenke's approach in
\cite{pom} is based on a slightly different, but equivalent,
renormalization which replaces \eqref{valiron-method}.
The approach in \cite{pom2}, which works for random iteration
sequences, needs some regularity hypothesis. On the other hand,
Cowen's construction \cite{Co} is based on an abstract model
relying strongly on the Riemmann uniformization theorem.
Finally, Bourdon and Shapiro's construction is based upon a
different renormalization process which works only with
some further regularity of $\v$ at $\tau$, but also guarantees
some stronger regularity properties for the semi-conjugation $\theta$.

In \cite[Prop. 6]{Br-Po} the first and last named authors
proved that actually all those different methods (when
applicable) provide essentially the same solution. Namely, if
$\tilde{\sigma}:\D\to\Hr$ is another (nonconstant) solution of
the functional equation \eqref{schroder} then there exists
$\lambda>0$ such that $\tilde{\sigma}=\lambda \sigma$.

Moreover, Valiron showed that $\sigma$ comes with some guaranteed, but
weak, regularity properties at $\tau$. In function theory language,
$\sigma$ is semi-conformal (or isogonal) at $\tau$, namely,
$\sigma$ fixes $\infty\in \de\Hr$ non-tangentially and
$\Arg\sigma$ has non-tangential limit $0$ at $\infty$.
As showed in \cite{Br-Po}, the semi-conformality of $\sigma$ is
essentially responsible for the uniqueness properties of $\sigma$ and
for the following dynamical properties of
$\phi$: for every orbit $z_n\defeq \phi^{\circ n}(z_0)$, $\Arg z_n$ tends to a
limit $\alpha(z_0)\in (-\pi/2,\pi/2)$ which depends harmonically on
$z_0$, and conversely, given an angle $\alpha\in (-\pi/2,\pi/2)$ one
can always find an orbit whose limiting argument is $\alpha$.

\subsection{Valiron's method in higher dimensions}

In $\C^N$, $N=2,3,...$, we let $\pi_j:\C^N\rightarrow\C$,
$j=1,...,N$, be the coordinate mappings; the usual inner product is
$\la z_1, z_2 \ra\defeq\sum_{j=1}^N z_{1,j} \overline{z_{2,j}}$, where
$z_{n,j}=\pi_j(z_n)$; the norm is $\|z\|^2\defeq\la z,z\ra$. The unit
ball $\B^N$ is $\{z\in \C^N: \|z\|^2<1\}$.

Let $\v$ be a holomorphic self-map of $\B^N$.
If $\v$ has no fixed points in $\B^N$ then B. MacCluer
\cite{Mac} proved that the Denjoy-Wolff theorem still holds.
Namely, the sequence of iterates of $\v$, $\{\v^{\circ n}\}$, converges
uniformly on compacta to the constant map $z\mapsto \tau$, $\forall
z\in \B^N$, for a
(unique) point $\tau\in\partial{\B^N}$ (called again the {\sl
Denjoy-Wolff point} of $\v$). Like in the one-dimensional
case, the number
\[
c :=\liminf_{z\rightarrow \tau}\frac{1-\|\v (z)\|}{1-\|z\|},
\]
belongs to $(0,1]$ and is called the {\sl multiplier} of $\v$
or the {\sl boundary dilatation coefficient} of $\v$ at $\tau$.
Also, $\tau$ is a fixed point in the sense of non-tangential
limits (and actually in the sense of $K$-limits as we define
below). However, in this case the differential of $\v$ might
not have nontangential limit at $\tau$. The map $\v$ is called
{\sl hyperbolic} if $c<1$ and {\sl parabolic} if $c=1$.

Here too $\v$ preserves certain ellipsoids internally tangent to $\de \B^N$
at $\tau$: defining
\begin{equation}\label{eq:horoball}
E(t)\defeq\left\{z\in\B^N: \frac{|1-\la z,
 \tau\ra|^2}{1-\|z\|^2}<1/t\right\},\end{equation}
then $\v(E(t))\subset E(t/c)$. In formulas,
\begin{equation}\label{eq:julia}
\frac{|1-\la \v(z),\tau\ra|^2}{1-\|\v(z)\|^2}\leq c\frac{|1-\la
  z,\tau\ra|^2}{1-\|z\|^2},
\end{equation}
for every $z\in \B^N$.

Assuming some regularity for $\v$ at $\tau$, in the spirit of
Bourdon-Shapiro, in \cite{Br-Ge} the first and the second named
authors proved that, if $\v$ is hyperbolic, one can solve the
following functional equation:
\[
\sigma \circ \v=A \sigma,
\]
where $\sigma:\B^N\to\C^N$ is a nonconstant holomorphic map with good
regularity properties at $\tau$, and where
$A$ is the matrix $d\v_\tau$. Recently such a result has been improved in
$\B^2$ by F. Bayart assuming less regularity for $\v$ at $\tau$
(see \cite{Ba} where also the parabolic case is considered).

On the other hand, the first and third named author in
\cite{Br-Po} have shown, that for all hyperbolic self-maps
$\v$, i.e., with no regularity assumptions at $\tau$, and for
each orbit $z_n=\v^{\circ n}(z_0)$, there is a Koranyi region
$K(\tau,R)$ such that $z_n$ will tend to $\tau$ while staying
in $K(\tau,R)$. Recall that, for $R>1/2$,  the {\sl $R$-Koranyi
approach region} at $\tau$ is a region of the form
\begin{equation}\label{koranji}
K(\tau,R)\defeq\{z\in\B^N: |1-\la z,\tau\ra|<R (1-\|z\|^2)\}.
\end{equation}

The original aim, when looking for semi-conjugations in the one-dimensional
case, was to show that general hyperbolic self-maps do indeed have a
similar dynamical behavior as the hyperbolic
automorphisms that share the same attracting fixed point.

In higher dimensions however, it is easy to construct maps whose image
lies in a sub-variety with non-zero codimension, and thus automorphisms alone
don't seem to be enough to model the dynamics of such maps (although
one may try to consider automorphisms of lower dimensional balls).
Also the fact that the differential of $\v$ does not in general have a
 non-tangential limit at $\tau$, shows that before trying to
 semi-conjugate $\v$ to an automorphism on an higher-dimensional ball,
 it is preferable to study
 the following ``one-dimensional'' equation first.
\begin{problem}\label{prob:onedim}
Find a nonconstant holomorphic map
$\Theta:\B^N\to \Hr\subset\C$ such that
\begin{equation}\label{valiron-piu-dim}
\Theta\circ \v=\frac{1}{c} \Theta.
\end{equation}
\end{problem}

The aim of this paper is to try to solve Problem \ref{prob:onedim} by
generalizing the method of Valiron to higher dimensions.

As in the one-dimensional
case, it is more convenient to move to the Siegel domain
\begin{equation}\label{eq:siegel}
\Hr^N:=\{(z,w)\in\C\times\C^{N-1}: \Re z>\|w\|^2\}
\end{equation}
which is biholomorphic to $\B^N$ via the Cayley transform
$\mathcal C:\B^N\to \Hr^N$ defined as
\begin{equation}\label{eq:caley}
\mathcal C(\zeta_1, \zeta'):=\left(\frac{1+\zeta_1}{1-\zeta_1},
\frac{\zeta'}{1-\zeta_1}\right).
\end{equation}
Thus, if
$\phi:\Hr^N\to \Hr^N$ is a hyperbolic holomorphic map with
Denjoy-Wolff point $\infty$ and multiplier $1/c$, we
define the following sequence
\begin{equation}\label{valivali}
\sigma_n(z,w):=\frac{\pi_1\circ \phi^{\circ n}(z,w)}{x_n},
\end{equation}
where,  $\pi_1(z,w):=z$ is the projection on the first
component and $x_n=\Re \pi_1 (\phi^{\circ n}(1,0))$. For short we will
say that the {\sl Valiron method works} whenever the sequence
$\{\sigma_n\}$ converges uniformly on compacta.

Our main result is the following:

\begin{maintheorem}
Let $\v:\B^N\to\B^N$ be a hyperbolic holomorphic self-map with
Denjoy-Wolff point $\tau\in\de\B^N$ and multiplier $c<1$. If
\begin{enumerate}
  \item there exists $z_0\in\B^n$ such that the sequence
  $\{\v^{\circ n}(z_0)\}$ is {\sl special} and
  \item the $\displaystyle{\hbox{K-}\lim_{z\to \tau} \frac{1-\la \v(z),\tau\ra}{1-\la z, \tau\ra}}$ exists,
\end{enumerate}
then the Valiron method works and there exists a nonconstant
holomorphic function $\Theta:\B^N\to \Hr$ such that
$\Theta\circ \v=\frac{1}{c}\Theta$.
\end{maintheorem}

In order to explain our hypotheses (1) and (2), we recall that
a sequence $\{z_n\}\subset \B^N$ converging to a point
$\tau\in\de\B^N$ is said to be {\sl special} if
\[
\lim_{n\to \infty}\frac{\|z_n-\la z_n,\tau\ra \tau\|^2}{1-|\la
z_n,\tau\ra|^2}=0,
\]
or, equivalently,  the Kobayashi distance $k_{\B^N}(z_n, \la
z_n,\tau\ra\tau)$, between $\{z_n\}$ and the projection of
$z_n$ along $\tau$, tends to zero as $n\to \infty$.
For the definition and properties of the Kobayashi distance we
refer to \cite{Kob} or \cite{Abate}; we will only use the
fact that the Kobayashi distance is invariant under biholomorphisms
and that $k_{\B^N}(0,z)=\tanh^{-1}(\|z\|)$.

Moreover,
a function $h:\B^N\to \C$ has {\sl K-limit} $L$ at
$\tau\in\de\B^N$, K-$\lim_{z\to \tau} h(z)=L$, if for any $R>1/2$
and any sequence $\{z_n\}\subset K(\tau, R)$ converging to
$\tau$ it follows that $\lim_{n\to \infty}h(z_n)=L$ (see
\cite{Abate} or \cite{Ru}).

Notice that if $\v:\B^N\to\B^N$ is a hyperbolic holomorphic
self-map with Denjoy-Wolff point $\tau\in\de\B^N$ and
multiplier $c<1$, then  Rudin's version of the classical
Julia-Wolff-Caratheodory theorem (see \cite[Thm. 8.5.6]{Ru} or
\cite[Thm. 2.2.29]{Abate}) implies that
\begin{equation}\label{eq:rjwc}
\lim_{n\to\infty}\frac{1-\la \v(z_n),\tau\ra}{1-\la z_n,
\tau\ra}=c
\end{equation}
for all sequences $\{z_n\}\subset\B^N$ converging to $\tau$
such that $\{z_n\}$ is special and $\{\la z_n, \tau\ra\}$
converges to $1$ nontangentially in $\D$. Such a limit is
called {\sl restricted K-limit}. Unfortunately, it is easy to show
that K-limits imply restricted K-limits, but not the converse. Thus,
hypothesis (2) is a non-trivial requirement.

Condition (1) is not always easy to verify, unless, say, the
map $\v$ happens to fix (as a set) a slice ending at $\tau$.
For instance, under the regularity assumptions of \cite{Br-Ge}
it follows that (2) holds, but it is not clear, {\em ex ante},
that (1) must also hold. On the other hand, once the
semi-conjugation is established in \cite{Br-Ge}, with good
regularity properties, then it is easy to verify that (1) had
to hold, {\em ex post}. In fact, we don't know of any explicit
examples where (1) fails. So it could be the case that (1) is
actually a superfluous hypothesis for the Main Theorem.

\subsection{An example}\label{ssec:example}
The following is an example of a map as in
the Main Theorem satisfying condition (1) but not (2)
and for which the Valiron method still works.

Consider the map
\[
\phi:\Hr^{2}\ni(z,w)\mapsto (Az+Aw^2\psi(z), 0)
\]
where $\psi:\Hr\to \D$ is any holomorphic function and $A>1$.
Then clearly, $\phi(\Hr^2)\subset \Hr^2$, $\infty$ is the
Denjoy-Wolff point of $\phi$, the multiplier is $A>1$  and the
sequence $\{\phi^{\circ n}(1,0)\}=\{(A^n,0)\}$ is special.
Moreover,
\[
\phi^{\circ n}(z,w)=A^n z+A^n w^2 \psi(z).
\]
Hence
\[
\sigma_n(z,w):=\frac{\pi_1\circ \phi^{\circ n}(z,w)}{x_n}=z+w^2\psi(z).
\]
Therefore $\{\sigma_n\}$ does not depend on $n$ and it can be
checked that the map $\sigma(z,w):=z+w^2\psi(z)$ solves $\sigma
\circ \phi=A\sigma$. Thus the Valiron method works. However, the
$K$-limit of $\frac{\phi_1(z,w)}{z}$ at $\infty$ does not exist
if $\psi$ doesn't have a non-tangential limit at $\infty$. In
particular, for such $\psi$, hypothesis (2) in the Main Theorem
is not satisfied.

It is interesting to note that for such an example, the crucial
equation \eqref{equi1} below becomes
\[
\frac{Ax_n z+Ax_n w^2\psi(x_nz)}{x_nz}=A+A\frac{w^2}{z}\psi(x_n
z),
\]
and the limit for $n\to \infty$ does not exists if $w\neq 0$.

In particular, the regularity hypothesis (2) in the Main
Theorem, while necessary in our proof, is not necessary for
Valiron's method to work.

Our Main Theorem is proved in Section \ref{sec:proof}.

In order to prove it, in Section \ref{limiti} we introduce a
new characterization of K-limits for functions, which we then
develop in the Appendix into the notion of {\sl E-limits}. We
believe that the new understanding of K-limits which comes from
the study of our E-limits might be a useful tool for other
results. In the last section we include some further comments
and open questions.

\section{Preliminaries on K-Limits}\label{limiti}

As mentioned before, we work in the Siegel domain (\ref{eq:siegel}).
A direct computation using \eqref{koranji} and (\ref{eq:caley}) shows that the Koranyi
region $K(\tau,R)$ with vertex at $\tau$ and amplitude $R$ in $\B^N$
corresponds to one with vertex at $\infty$ and amplitude $M\defeq 2R>1$ in $\Hr^N$
given by
\begin{equation}\label{eq:hkoranyi}
K(\infty,M)\defeq\left\{(z,w)\in\Hr^N: \|w\|^2 < \Re z - \frac{|z+1|}{M}\right\}.
\end{equation}
To get a geometric feeling for these objects, notice that the
ellipsoids $E(t)$ defined in (\ref{eq:horoball}) correspond in $\Hr^N$
to the sets
\[
\Eo(T)\defeq\{(z,w)\in \Hr^N: \Re z -\|w\|^2 > T\}
\]
for some $T>0$ fixed. So, in particular, a sequence in $K(\infty,M)$
tending to infinity will eventually be contained in every $\Eo(T)$
for $T$ large, because $z$ tends to infinity when $(z,w)\in \Hr^{N}$
tends to infinity.

Notice also that the property (\ref{eq:julia}) for a hyperbolic map
$\phi:\Hr^N\rightarrow\Hr^N$ with multiplier $A>1$ reads as follows:
\[
\Re \phi_1(z,w) -\|\phi^\prime(z,w)\|^2 > A(\Re z-\|w\|^2)
\]
for every $(z,w)\in \Hr^N$.

We will find it convenient to use an equivalent characterization of
K-limits. First we need a few definitions.

For $Z=(z,w)\in \Hr^N$, let $p(Z)\defeq (z,0)$ be the projection of
$Z$ onto the complex line $L\defeq\{(z,0): z\in \Hr\}\subset\Hr^N$

\begin{definition}
Let $Z_{n}= (z_{n},w_{n})\in\Hr^N$ converge to $\infty$.
\begin{itemize}
\item[(i)] We say the
convergence is {\sl $C$-special} if there exists $0\leq C <\infty$ such that
\[
k_{\Hr^N}(Z_n, p(Z_n))\leq C, \qquad \forall n,
\]
where $k_{\Hr^{N}}$ is the Kobayashi distance on $\Hr^{N}$.
\item[(ii)]
We say the convergence is {\sl restricted} if $\{z_{n}\}$ converges
non-tangentially to $\infty$ in $\Hr$.
\end{itemize}
\end{definition}
\begin{remark}
The concepts just introduced of $C$-special and restricted
sequences are formulated using the complex geodesic $z\in
\Hr\mapsto (z,0)\in\Hr^{N}$ and the projection associated to
it. It turns out that being $C$-special and restricted
 do not depend on the chosen complex geodesic  with
$\infty$ in its boundary. This is used in the proof of the Main
Theorem and could be useful in domains other than $\Hr^{N}$ and
$\B^N$. For this reason, in the Appendix, Section
\ref{sec:appendix}, we provide a rigorous proof of this fact.
\end{remark}

\begin{remark}
A $0$-special sequence is simply referred to as {\sl special},
see also \cite{Abate} and \cite{Ru}.
\end{remark}

\begin{lemma}\label{lem:klimits}
Let $Z_{n}= (z_{n},w_{n})\in\Hr^N$ converge to $\infty$. Then, the
following are equivalent:
\begin{enumerate}
\item $Z_{n}$ stays inside a Koranyi region $K (\infty ,M)$ for some
  $1<M<\infty$;
\item $Z_{n}$ is $C$-special, for some $C<\infty$, and is restricted;
\item There is $0<a<1$ and $0<T<\infty$, such that
\[
\|w\|^2\leq a\Re z_n\qquad\mbox{ and }\qquad |\Im z_n| \leq T \Re z_n.
\]
\end{enumerate}
\end{lemma}
The proof of Lemma \ref{lem:klimits} rests on the following
computation. For $Z=(z,w)\in\Hr^N$, we compute  the Kobayashi distance in
$\Hr^N$ between $Z$ and $p(Z)$. Set $z=x+iy$ and
notice that the map $T(u,v)=(\frac{u-iy}{x},\frac{v}{\sqrt{x}})$ is
an automorphism of $\Hr^N$. Thus by invariance, we have
\begin{equation}\label{kob-iperplane}
\begin{split}
    k_{\Hr^N}((z,0),(z,w))&=k_{\Hr^N}((1,0),
    T(z,w))=k_{\B^N}(0, \mathcal
    C^{-1}(T(z,w)))\\&=\tanh^{-1}
    \|\mathcal
    C^{-1}(T(z,w))\|=\tanh^{-1}\|(0,\frac{w}{\sqrt{x}})\|
    =\tanh^{-1}\frac{\|w\|}{\sqrt{x}}.
\end{split}
\end{equation}
In other words, $k_{\Hr^N}(Z,p(Z))=\tanh^{-1}(\|w\|/\sqrt{\Re z})$ and
it is useful to
recall that $\tanh^{-1}(s)=(e^s-e^{-s})/(e^s+e^{-s})$ is a positive increasing function on
$(0,1)$ with a vertical asymptote at $1$.
\begin{proof}[Proof of Lemma \ref{lem:klimits}]
By \eqref{kob-iperplane}, a sequence $Z_n=(z_n,w_n)\in \Hr^N$ is $C$-special
for some $0<C<\infty$ if and only if
$\|w_n\|^2 \leq a \Re z_n$ for some $0<a<1$. In fact, $a=\tanh C$.
Thus, since $\Im z_n \leq T\Re z_n$ is an usual formulation of
non-tangentiality in $\Hr$, we have that (2) and (3) are
equivalent.

Assuming (3) and writing $z_n=x_n+iy_n$, we have $|z_n+1|^2\leq
(1+T^2)x_n^2+2x_n+1$. Thus
\[
x_n-\frac{|z_n+1|}{M}\geq \left(1-\frac{\sqrt{1+T^2}}{M}\right)x_n+o\left(\frac{1}{x_n}\right),
\]
as $x_n$ tends to infinity. Choose $M$ large enough, so that
$1-\sqrt{1+T^2}/M<a<1$. This
ensures that $Z_n\in K(\infty, M)$ for all $n$ large. So (3) implies (1).

Conversely, assume that $Z_n\in K(\infty,M)$ for some
$1<M<\infty$. Then, since
\[
x_n-|z_n+1|/M\leq (1-1/M)x_n,
\]
by (\ref{eq:hkoranyi}), we have $\|w_n\|^2\leq
a\Re z_n$ with $a=1-1/M$. Also, $|z_n+1|\leq M\Re z_n$, so $|\Im
z_n|\leq M\Re z_n$. Hence, (1) implies (3).
\end{proof}

\section{The proof of the Main Theorem}\label{sec:proof}

We start by reformulating it in the context of $\Hr^{N}$.

\begin{maintheorem}[Siegel domain version]
Let $\phi=(\phi_1,\phi'):\Hr^N\to\Hr^N$ be holomorphic, with
Denjoy-Wolff point $\infty$ and multiplier $\lambda>1$. Assume
that
\begin{enumerate}
    \item There exists $Z_0\in\Hr^N$ such that the sequence
    $\{\phi^{\circ n}(Z_0)\}$ is special.
    \item  $\hbox{K-$\lim$}_{\Hr^N\ni (z,w)\to\infty}
    \frac{\phi_1(z,w)}{z}$ exists.
\end{enumerate}
Then Valiron's method works and there exists a non-constant
holomorphic map $\sigma:\Hr^N\to\Hr$ such that
\[
\sigma \circ \phi=\lambda \sigma.
\]
\end{maintheorem}
\begin{remark}\label{rem:conj}
By considering $T\circ \phi \circ T^{-1}$, where $T$ is an
automorphism of $\Hr^N$ fixing $\infty$ and such that $T (Z_{0})=
(1,0)$, we can always assume that it is the sequence $\phi^{\circ n}
(1,0)$ that is special, see the proof of Lemma
\ref{specialspecial}. So we will make this assumption in the sequel.
\end{remark}

\begin{remark}\label{rem:valeVali}
The Valiron method is invariant under conjugation, namely, let
$\phi:\Hr^N\to\Hr^N$ be hyperbolic holomorphic with
Denjoy-Wolff point $\infty$, let $T$ be an automorphism of
$\Hr^N$ fixing $\infty$ and let $\tp:=T\circ \phi\circ T^{-1}$.
Then the sequence $\{\sigma_n\}:=\{(\pi_1\circ \phi^{\circ
n})/x_n\}$ given by \eqref{valivali} converges  if and only if
the sequence $\{\ts_n\}:=\{(\pi_1\circ \tp^{\circ
n})/\tilde{x}_n\}$ converges (here $\tilde{x}_n=\Re
\pi_1(\tp^{\circ n}(1,0))=\Re \pi_1(T(\phi^{\circ n}(
T^{-1}(1,0))))$). In fact, by a direct computation, it turns
out that if $\sigma_n\to \sigma$ as $n\to \infty$ then
$\ts_n\to (x_0-\|w_0\|^2)\sigma \circ T^{-1}$, where
$(x_0+iy_0, w_0):=T^{-1}(1,0)$. We leave the details of such a
computation to the reader.
\end{remark}

We need a preliminary result.

\begin{lemma}\label{real and image}
Let $\phi=(\phi_1,\phi'):\Hr^N\to\Hr^N$ be holomorphic, with
Denjoy-Wolff point $\infty$ and multiplier $\lambda\geq 1$.
Assume the sequence $\{\phi^{\circ n}(1,0)\}$ is special. Write
$\phi^{\circ n}(1,0)=(z_n,w_n)$ and $z_{n}=x_{n}+iy_{n}$. Then
\begin{enumerate}
\item $\displaystyle{\lim_{n\to\infty}\frac{x_{n+1}}{x_n}=\lambda}$.
\item There exists $L\in \R$ such that $\displaystyle{ \lim_{n\to \infty}\frac{y_n}{x_n}=L.}$
\end{enumerate}
\end{lemma}
\begin{proof}
As proved in \cite[section
3.5]{Br-Po}, for any
fixed $Z\in\Hr^N$, the orbit $\{\phi^{\circ n}(Z)\}$ stays in a
Koranyi region with vertex at $\infty$ and so, in particular, it is
restricted. Therefore, there exists
$C>0$ such that for all $n\in\N$
\begin{equation}\label{restricted}
|y_n|\leq C x_n.
\end{equation}
By Rudin's version of the classical
Julia-Wolff-Caratheodory theorem (\ref{eq:rjwc}), reformulated in
$\Hr^{N}$ (see Theorem \ref{JWC} in the Appendix), since
$(z_{n},w_{n})$ is special and restricted, it follows that
\[
\lim_{n\to
\infty}\frac{z_{n+1}}{z_n}=\frac{\phi_1(z_n,w_n)}{z_n}=\lambda.
\]
In particular we can write
\begin{equation}\label{znuno}
z_{n+1}=\lambda z_n+ o(1)z_n.
\end{equation}
 Dividing \eqref{znuno} by $x_n$
and taking the real part,  we obtain $\frac{x_{n+1}}{x_n}=\lambda +
\Re o(1) - \frac{y_n}{x_n}\Im o(1)$. Taking the limit for $n\to
\infty$, by \eqref{restricted}, we get
\begin{equation}\label{limitxn}
   \lim_{n\to\infty}\frac{x_{n+1}}{x_n}=\lambda,
\end{equation}
which proves (1).

In order to prove (2), let $\left\{\frac{y_{n_k}}{x_{n_k}}\right\}$ be any
convergent subsequence and let $L$ be its limit. By
\eqref{restricted}, $L$ is finite. Moreover,
\begin{equation}\label{samelimit}
\frac{z_{n+1}}{z_n}=\frac{x_{n+1}}{x_n}\frac{1+i\frac{y_{n+1}}{x_{n+1}}}{1+i\frac{y_{n}}{x_{n}}}
\end{equation}
and by \eqref{znuno} and  \eqref{limitxn} we see that
$\left\{\frac{y_{n_k+1}}{x_{n_k+1}}\right\}$ is also
a convergent sequence with the
same limit $L$. Assume by contradiction that there exists a
converging subsequence $\left\{\frac{y_{m_k}}{x_{m_k}}\right\}$   with limit
$L'\neq L$. Let
\[
q_n:=\frac{x_{n+1}}{x_n}+i\frac{y_{n+1}-y_n}{x_n}.
\]
By \eqref{limitxn}, we have
\[
\Im
q_{n_k}=\frac{y_{n_k+1}-y_{n_k}}{x_{n_k}}=\frac{y_{n_k+1}}{x_{n_k+1}}\frac{x_{n_k+1}}{x_{n_k}}
-\frac{y_{n_k}}{x_{n_k}}\longrightarrow L(\lambda-1),
\]
and similarly $\Im q_{m_k}\to L'(\lambda-1)$. Therefore
$\{q_{n_k}\}$ converges to $\lambda+i L(\lambda-1)$ while
$\{q_{m_k}\}$ converges to $\lambda+i L'(\lambda-1)$.

We claim that $\{q_n\}$ can have at most two accumulation
points, say $a, a'$ (which must be necessarily
$a=\lambda+iL(\lambda-1)$ and $a'=\lambda+iL'(\lambda-1)$).
Assuming the claim is true, let $U, U'$ be two open
neighborhoods of $a$ and $a'$ respectively such that $U\cap
U'=\emptyset$. Since $\{q_n\}$ has only $a,a'$ as accumulation
points by our claim, there exists $n_0$ such that for all
$n>n_0$ then either $q_n\in U$ or $q_n\in U'$. Moreover, since
$\{q_{n_k}\}\subset U$ for $n_k>n_0$ and  $\{q_{m_k}\}\subset
U'$ for $m_k>n_0$, one can select a subsequence
$\{q_{l_k}\}\subset U$ such that $\{q_{l_k+1}\}\subset U'$. But
this implies that $\left\{\frac{y_{l_k}}{x_{l_k}}\right\}$ converges to
$L(\lambda-1)$ while $\left\{\frac{y_{l_k+1}}{x_{l_k+1}}\right\}$
converges to $L'(\lambda-1)$, contradicting our previous
argument in \eqref{samelimit}.

We are left to show that $\{q_n\}$ can have at most two accumulation
points. We already know that $\Re q_n\to \lambda>1$. We are going to
show that the (real) sequence $\{k_{\Hr}(1,q_n)\}$ of hyperbolic
distances between $1$ and $q_n$ has limit, say $d$. Thus the
accumulation points of $\{q_n\}$ must belong to the intersection
between the real line $\{\zeta\in \Hr: \Re \zeta=\lambda\}$ and the
boundary of the hyperbolic disc of center $1$ and radius $d$, and this
intersection
consists of at most two points.

To see that $\{k_{\Hr}(1,q_n)\}$ converges, let us introduce the
family of automorphisms of $\Hr^N$ given by
\begin{equation}\label{Tn}
T_n(z,w):=\left(\frac{z-iy_n}{x_n}, \frac{w}{\sqrt{x_n}}\right).
\end{equation}
Notice that $T_n(z_n,0)=(1,0)$ and $T_n\circ
T_{n+1}^{-1}(1,0)=(q_n,0)$, from which we obtain that
\begin{equation}\label{ineq0}
\begin{split}
k_{\Hr}(1,q_n)&=k_{\Hr^N}((1,0), (q_n,0))=k_{\Hr^N}((1,0),T_n\circ
T_{n+1}^{-1}(1,0))\\&=k_{\Hr^N}(T_n^{-1}(1,0),
T_{n+1}^{-1}(1,0))=k_{\Hr^N}((z_n,0),(z_{n+1},0)).
\end{split}
\end{equation}
Now, by the contracting property of Kobayashi's distance,
\begin{equation}\label{ineq1}
\begin{split}
k_{\Hr^N}((z_n,0)&,(z_{n+1},0))=k_{\Hr}(z_n,z_{n+1})\\&=k_{\Hr}(\pi_1(z_n,w_n),\pi_1(z_{n+1},
w_{n+1}))\leq k_{\Hr^N}((z_n,w_n),(z_{n+1},w_{n+1})).
\end{split}
\end{equation}
On the other hand, by the triangle inequality,
\begin{equation}\label{ineq2}
\begin{split}
k_{\Hr^N}((z_n,0),(z_{n+1},0))&\geq
k_{\Hr^N}((z_n,w_n),(z_{n+1},w_{n+1}))\\&-k_{\Hr^N}((z_n,0),(z_{n},w_{n}))-
k_{\Hr^N}((z_{n+1},0),(z_{n+1},w_{n+1})).
\end{split}
\end{equation}
Since $\{(z_n, w_n)\}$ is special then both
$k_{\Hr^N}((z_n,0),(z_{n},w_{n}))$ and
$k_{\Hr^N}((z_{n+1},0),(z_{n+1},w_{n+1}))$ tend to $0$ as $n\to
\infty$. Therefore, from \eqref{ineq0}, \eqref{ineq1} and
\eqref{ineq2} it follows
\[
\lim_{n\to \infty} k_\Hr(1,q_n)=\lim_{n\to
\infty}k_{\Hr^N}((z_n,0),(z_{n+1},0))=\lim_{n\to \infty}
k_{\Hr^N}((z_n,w_n),(z_{n+1},w_{n+1})),
\]
and the latter limit exists  because the sequence
$\{k_{\Hr^N}((z_n,w_n),(z_{n+1},w_{n+1}))\}$ is non-increasing
in $n$ since the Kobayashi distance is contracted by
holomorphic maps.
\end{proof}

\begin{proof}[Proof of the Main Theorem]
As mentioned in Remark \ref{rem:conj} and Remark
\ref{rem:valeVali}, after conjugating $\phi$ with some
automorphism of $\Hr^N$  we can suppose that $Z_0=(1,0)$, and
as we saw in the proof of Lemma \ref{real and image}, the orbit
of $(1,0)$ is thus both special and restricted. Moreover, see
Proposition \ref{tech2} in the Appendix, the conjugation made
does not effect our regularity hypothesis, namely
\begin{equation}\label{regularity}
{K\hbox{-}\lim}_{\Hr^N\ni
    (z,w)\to\infty}\frac{\phi_1(z,w)}{z}=\lambda.
\end{equation}
Letting
$(z_n,w_n):=\phi^{\circ n}(1,0),\ z_n=x_n+iy_n$,
and using Lemma \ref{lem:klimits}, we see that
\begin{equation}\label{speciale}
    \lim_{n\to\infty} \frac{\|w_n\|}{\sqrt{x_n}}=0.
\end{equation}

Now we consider the Valiron-like sequence $\{\sigma_n\}$ of
holomorphic maps from $\Hr^N$ to $\Hr$ defined by
\[
\sigma_n(z,w):=\frac{\pi_1\circ \phi^{\circ n}(z,w)}{x_n},
\]
where, as usual, $\pi_1(z,w):=z$ is the projection on the first
component. Notice that
\begin{equation}\label{intert}
\sigma_n \circ \phi = \frac{ \pi_1\circ
\phi^{\circ (n+1)}}{x_n}=\frac{x_{n+1}}{x_n} \sigma_{n+1}.
\end{equation}
If we can prove that the sequence $\{\sigma_n\}$ converges
uniformly on compacta  to a non-constant   map
$\sigma:\Hr^N\to\Hr$ (which is necessarily holomorphic), then
by  taking the limit for $n\to \infty$ in \eqref{intert}, and by
Lemma \ref{real and image}~(1), we obtain that $\sigma\circ
\phi=\lambda \sigma$.

We will now  show that $\{\sigma_n\}$ is uniformly convergent
on compacta to a non-constant function.

First of all, we notice that by Lemma \ref{real and image}~(2),
\[
\sigma_n(1,0)=1+i\frac{y_n}{x_n}\longrightarrow 1+i L,
\]
as $n\to\infty$. And, on the other hand, again by Lemma \ref{real
and image}
\[
\sigma_n(\phi(1,0))=\frac{\pi_1\circ
\phi^{\circ (n+1)}(1,0)}{x_n}=\frac{x_{n+1}+iy_{n+1}}{x_n}=
\frac{x_{n+1}}{x_n}+i\frac{x_{n+1}}{x_n}\frac{y_{n+1}}{x_{n+1}}\to
\lambda+i \lambda L,
\]
as $n\to \infty$. Since $\lambda>1$, the above proves that any
limit of the sequence $\{\sigma_n\}$ cannot be constant.

Now we are going to prove that for any $(z,w)\in\Hr^N$
\begin{equation}\label{specialsigma}
    \lim_{n\to\infty}k_\Hr(\sigma_n(z,w),\sigma_{n+1}(z,w))=0.
\end{equation}
To this aim, we first notice that the set $\{\sigma_n(z,w)\}$ is
relatively compact in $\Hr$. Indeed, let $\pi_w:\C^N\to \C^{N-1}$ be
the projection $\C\times \C^{N-1}\ni (z,w)\mapsto w\in \C^{N-1}$ and
define
\begin{equation}\label{Sn}
S_n(z,w):=\left(\sigma_n(z,w),
\frac{\pi_w(\phi^{\circ n}(z,w))}{\sqrt{x_n}}\right).
\end{equation}
Notice that $S_n=L_n \circ \phi^{\circ n}$, where $L_n$ is the
automorphism of $\Hr^N$ defined by $L_n(z,w)=(z/x_n,
w/\sqrt{x_n})$. Therefore $S_n:\Hr^N\to \Hr^N$. Moreover, by
Lemma \ref{real and image}~(1) and \eqref{speciale}
\begin{equation}\label{limitS}
S_n(1,0)=\left(\sigma_n(1,0), \frac{w_n}{\sqrt{x_n}}
\right)=(1+i\frac{y_n}{x_n}, \frac{w_n}{\sqrt{x_n}})\to (1+iL,0),
\end{equation}
as $n\to \infty$.  In particular there exists $C>0$ such that
$k_{\Hr^N}(S_n(1,0), (1+iL,0))<C$ for all $n\in \N$. Therefore, by the
triangle inequality and the contraction property,
\begin{equation*}
\begin{split}
k_{\Hr^N}(S_n(z,w), (1+iL,0))&\leq k_{\Hr^N}(S_n(z,w),
S_n(1,0))+k_{\Hr^N}(S_n(1,0), (1+iL,0))\\&\leq
k_{\Hr^N}((z,w),(1,0))+C,
\end{split}
\end{equation*}
which proves that $\{S_n(z,w)\}$ is relatively compact in $\Hr^N$.

Now,  notice that
\[
\sigma_{n+1}=\pi_1 \circ L_{n+1}\circ \phi \circ L_n^{-1}\circ S_n.
\]
Since we already proved that the sequence $\{S_n(z,w)\}$ is
relatively compact in $\Hr^N$, \eqref{specialsigma} will follow if
we prove that $\pi_1 \circ L_{n+1}\circ \phi \circ L_n^{-1}\to \pi_{1}$
as $n\to \infty$. A direct computation shows that
\begin{equation}\label{equi1}
\pi_1 \circ L_{n+1}\circ \phi \circ
L_n^{-1}(z,w)=\frac{\pi_1(\phi(x_nz,
\sqrt{x_n}w))}{x_{n}z}\frac{x_nz}{x_{n+1}}.
\end{equation}
Now for all $n\in \N$, by \eqref{kob-iperplane}
\[
k_{\Hr^N}((x_nz, \sqrt{x_n}w),
(x_nz,0))=\tanh^{-1}\frac{\|w\|}{\sqrt{\Re z}}<\infty,
\]
and clearly $\{x_nz\}$ converges to $\infty$ non-tangentially in $\Hr$.
Thus the sequence $\{(x_nz, \sqrt{x_n}w)\}$ is $C$-special and
restricted. Hence, by applying \eqref{regularity} and Lemma
\ref{real and image}~(1) to the limit as $n\to \infty$ in
\eqref{equi1}, we get $\pi_1 \circ L_{n+1}\circ \phi \circ
L_n^{-1}(z,w)\to z$ as $n\to \infty$, as needed.

At this point, let $\{\sigma_{n_k}\}$ be a convergent
subsequence of $\{\sigma_n\}$ and let $\sigma$ be its limit, which we
know is non-constant. By
\eqref{specialsigma}, $\{\sigma_{n_k+1}\}$ also converges to
$\sigma$. By \eqref{intert} and Lemma \ref{real and image}~(1)
we see that
\begin{equation}\label{functeq}
\sigma \circ \phi = \lambda \sigma.
\end{equation}
It remains to show that the Valiron method works, namely, that
the sequence $\{\sigma_n\}$ converges. By  the very definition,
$\{\sigma_n\}$ converges if and only if $\{\pi_1\circ S_n\}$
does---with $S_n$ defined in \eqref{Sn}. We already proved that
$\{S_n\}$ is bounded on compacta of $\Hr^N$, thus it is a
normal family. Let $S$ be a limit of $\{S_n\}$. Let
$Z\in\Hr^N$. Since the Kobayashi distance is contracted by
holomorphic maps, the sequence $\{k_\Hr(S_n(1,0),S_n(Z))\}$ is
decreasing in $n$ and must have a limit. Therefore, by
\eqref{limitS}, for all $Z\in \Hr^N$,
\[
\lim_{n\to \infty}k_\Hr(S_n(1,0),S_n(Z))=k_\Hr((1+iL,0), S(Z)).
\]
This implies that if $\tilde{S}$ is another limit of $\{S_n\}$
then $k_\Hr((1+iL,0), S(Z))=k_\Hr((1+iL,0), \tilde{S}(Z))$ for
all $Z\in \Hr^N$. Thus, conjugating both $S, \tilde{S}$ with a
Cayley map $\mathcal C'$ which maps $(1+iL,0)$ into $O\in
\B^N$, we find two holomorphic maps $S', \tilde{S}':\B^N\to
\B^N$ with the property that $\|S'(Z)\|=\|\tilde{S}'(Z)\|$ for
all $Z\in \B^N$. Hence (see, {\sl e.g.}, \cite[Prop. 3 p.
102]{Da}) there exists a unitary matrix $U$ such that
$S'=U\tilde{S}'$. Translating into $\Hr^N$ this means that
$\tilde{S}=T\circ S$ for some automorphism $T:\Hr^N\to\Hr^N$
fixing $(1+iL,0)$. We claim that $\pi_1\circ T(z,w)=z$, hence
$\pi_1\circ S=\pi_1\circ \tilde{S}$ which implies that
$\{\pi_1\circ S_n\}$---and hence $\{\sigma_n\}$---is
converging. In order to prove that $\pi_1\circ T(z,w)=z$, it is
enough to prove that $T(z,0)=(z,0)$ for some point $z\in \Hr\setminus
\{1+iL\}$,
because then by the classical theory of automorphisms (see
\cite{Abate} or \cite{Ru}) $T$ must fix pointwise the complex
geodesic $\Hr\times\{0\}$. To this aim, let $Z_1:=\phi(1,0)$. Let
$\{S_{n_k}\}$ be a sub-sequence of $\{S_n\}$ converging to $S$.
By \eqref{functeq},
\[
(\pi_1\circ S)(Z_1)= (\pi_1\circ S)(\phi(1,0))=\lambda
\sigma(1,0)=\lambda (1+iL).
\]
On the other hand, setting as before $(z_n,w_n):=\phi^{\circ n}(1,0)$, we
get
\begin{equation*}
\begin{split}
(\pi_w \circ
S)(Z_1)&=\lim_{k\to\infty}\frac{\pi_w(\phi^{\circ n_k}(Z_1))}{\sqrt{x_{n_k}}}
=\lim_{k\to\infty}\frac{\pi_w(\phi^{\circ (n_k+1)}(1,0))}{\sqrt{x_{n_k}}}\\
&=\lim_{k\to\infty}\frac{w_{n_k+1}}{\sqrt{x_{n_k}}}=
\lim_{k\to\infty}\frac{w_{n_k+1}}{\sqrt{x_{n_k+1}}}\sqrt{\frac{x_{n_k+1}}{x_{n_k}}}=0,
\end{split}
\end{equation*}
where the last equality follows from \eqref{speciale} and Lemma
\ref{real and image}~(1). Thus $S(Z_1)=(\lambda (1+iL),0)$.
Similarly, we have $\tilde{S}(Z_1)=(\lambda (1+iL),0)$.
Therefore
\[
T(\lambda (1+iL),0)=(T\circ S)(Z_1)=\tilde{S}(Z_1)=(\lambda
(1+iL),0),
\]
which proves that $\pi_1\circ T(z,0)=z$ as needed.
\end{proof}

\section{Further remarks and open questions}

{\bf 1.} In order to make the Valiron construction to work, in
the Main Theorem we need the technical hypothesis (1),
namely  that $\phi$ possesses a $0$-special orbit. We do not know
whether any hyperbolic holomorphic self-map of the ball always
has such an orbit or not. Clearly, if the self-map has an
invariant complex geodesic (whose closure must necessarily
contain the Denjoy-Wolff point) then such a condition is
satisfied for all points on such a complex geodesic. For
instance, if $T:\B^N\to \B^N$ is a hyperbolic automorphism with
Denjoy-Wolff point $e_1$ and other fixed point  $-e_1$, then
the orbit of any point $(z,0')$ is (obviously) special, and
conversely, the orbit of any point of the form $(z,z')$ with
$z'\neq 0'$ is not special.

{\bf 2.} As shown by Example in section \ref{ssec:example},
hypothesis (2) in the Main Theorem  is not necessary for
the Valiron construction to work in higher dimension.

{\bf 3.} Along the lines of the one-dimensional Valiron
construction (see, {\sl e.g.}, \cite[p. 47]{Br-Po}) one can
prove that if $\sigma$ is the intertwining map given by the Main
Theorem, then  $\Hr\ni \zeta\mapsto \sigma(\zeta,0)$ is {\sl
semi-conformal} at $\infty$. However, no further regularity on
$\sigma$ at $\infty$ seems to follow from the construction.

{\bf 4.} Uniqueness (up to composition with linear fractional
maps) of intertwining mappings in higher dimension---without
assign further conditions---does not hold. The main theoretical
reason is that in dimension one the centralizer of a given
hyperbolic automorphism consists of hyperbolic automorphisms
while in higher dimension this is not longer so (see
\cite{Ge-de}). For example, if $H:\B^N\to \B^N$ is a hyperbolic
automorphism, then any holomorphic self-map $F:\B^N\to \B^N$
such that $F\circ H=H\circ F$  solves the (trivial) Schr\"oder
equation $\sigma \circ H=H\circ \sigma$. By \cite{Ge-de}, if
$N>1$, then there exist mappings $F$ which are not linear
fractional maps.

\section{Appendix: E-limits}\label{sec:appendix}

In this appendix, we introduce the notion of $E$-limit in $\Hr^{N}$
and show that is equivalent to that of $K$-limit. However, this new
definition might be useful in more general domains. We also prove a
couple of routine facts that were needed in the proof of the Main Theorem.

A {\sl complex geodesic} $f:\Hr\to \Hr^N$ is a holomorphic map
which is an isometry between the Poincar\'e distance on $\Hr$
and the Kobayashi distance on $\Hr^N$. It is well known (see,
{\sl e.g.}, \cite{Abate}) that for $\Hr^N$ the image  of a
complex geodesic is the intersection of $\Hr^N$ with an affine
complex line. A {\sl linear projection} $\rho:\Hr^N\to \Hr^N$
is a holomorphic map such that $\rho^2=\rho$, the image
$\rho(\Hr^N)$ is the intersection of $\Hr^N$ with an affine
complex line (namely it is a complex geodesic) and
$\rho^{-1}(\rho(Z))$ is an affine hyperplane in $\Hr^N$ for all
$Z\in\Hr^N$. To any complex geodesic it is associated a unique
linear projection and conversely, to any linear projection it
is associated a unique (up to parametrization) complex
geodesic.

Given any complex geodesic $f:\Hr\to \Hr^N$ there exists an
automorphism $G$ of $\Hr^N$ such that $f(\zeta)=G^{-1}(\zeta,
0)$. The  linear projection associated to $f$ is then given by
$\rho(z,w)=G^{-1}(\pi_1(G(z,w)),0)$, where $\pi_1(z,w):=z$. The
map $\tr:=f^{-1}\circ \rho:\Hr^N\to \Hr$ is called the {\sl
left inverse} of $f$.

If $\rho:\Hr^N\to \Hr^N$ is a linear projection such that
$\overline{\rho(\Hr^N)}$ contains $\infty$, for short we say
that $\rho$ is a linear projection at $\infty$.

We will denote by $p_1:\Hr^N\to \Hr^N$ the linear projection at
$\infty$ given by $p_1(z,w)=(z,0)$, associated to the complex
geodesic $f(\zeta)=(\zeta,0)$ and left inverse $\pi_1(z,w)=z$.

\begin{definition}
Let $\rho:\Hr^N\to \Hr^N$ be a linear projection at $\infty$. A
sequence $\{Z_k\}\subset\Hr^N$ converging to $\infty$ is said
{\sl $C$-special with respect to $\rho$} if there exists $C\geq
0$ such that
\[
\limsup_{k\to \infty} k_{\Hr^N}(Z_k, \rho(Z_k))\leq C.
\]
The sequence $\{Z_k\}$ converging to $\infty$ is said to be {\sl
$\rho$-restricted} if $\{\rho(Z_k)\}$ converges non-tangentially to
$\infty$ in $\rho(\Hr^N)$.
\end{definition}

\begin{lemma}\label{specialspecial}
Let $\{Z_k\}\subset \Hr^N$ be a sequence converging to
$\infty$. Let $\rho_0:\Hr^N\to \Hr^N$ be a linear projection at
$\infty$. Then $\{Z_k\}$ is $C$-special ($C\geq 0$) with
respect to $\rho_0$ if and only if it is $C$-special (same $C$)
with respect to any linear projection at $\infty$ $\rho$. The
sequence $\{Z_k\}$ is $\rho_0$-restricted if and only if it is
$\rho$-restricted with respect to any  linear projection at
$\infty$ $\rho$.
\end{lemma}

\begin{proof}
Let $T_0$ be an automorphism of $\Hr^N$ fixing $\infty$ such
that $\rho(z,w)=T_0^{-1}(p_1(T_0(z,w))$. Since $T_0$ is an
isometry for $k_{\Hr^N}$,  then $\{Z_k\}$ is $C$-special with
respect to $\rho_0$ ({\sl respectively} $\rho_0$-restricted) if
and only if $\{T_0(Z_k)\}$ is $C$-special with respect to $p_1$
({\sl respect.} $p_1$-restricted). Therefore it is enough to
prove that if $\{Z_k\}$ is $C$-special with respect to $p_1$
({\sl respectively} $p_1$-restricted) then it is $C$-special
with respect to any linear projection at $\infty$ $\rho$ ({\sl
respect.} $\rho$-restricted).

Given a linear projection at $\infty$ $\rho$, there exists
$a\in \C^{N-1}$ and an automorphism $T\in\Aut(\Hr^N)$ of the
type
\[
T(z,w)=(z+\|a\|^2+2\la w,a\ra, w+a)
\]
such that $\rho=T^{-1}\circ p_1\circ T$. A direct computation
shows that
\begin{equation}\label{rho}
\rho(z,w)=(z+2\|a\|^2+2\la w, a\ra, -a).
\end{equation}
Therefore, writing $Z_k=(z_k,w_k)=(x_k+iy_k, w_k)$ and, arguing
similarly to \eqref{kob-iperplane}, we obtain
\begin{equation*}
\begin{split}
k_{\Hr^N}(p_1(Z_k),\rho(Z_K))&=k_{\Hr^N}((z_k,0),
(z_k+2\|a\|^2+2\la w_k, a\ra, -a))\\&=k_{\Hr^N}\left((1,0),
(\frac{z_k+2\|a\|^2+2\la w_k, a\ra-iy_k}{x_k},
\frac{-a}{\sqrt{x_k}})\right)\\&=\tanh^{-1}\sqrt{\frac{|2\|a\|^2+2\la
w_k,
a\ra|^2+4x_k\|a\|^2}{x_k^2\left|2+\frac{2\|a\|^2}{x_k}+2\frac{\la
w_k, a\ra}{x_k}\right|^2}}.
\end{split}
\end{equation*}
The last term tends to $0$ as $x_k\to \infty$, which is the
case if $k\to \infty$ because $Z_k\to \infty$ and $x_k=\Re
z_k>\|w_k\|^2$. Thus
\begin{equation}\label{specikilo}
    \lim_{k\to\infty}k_{\Hr^N}(p_1(Z_k),\rho(Z_K))=0.
\end{equation}
Now, using the triangle inequality and \eqref{specikilo} we see that
if $\{Z_k\}$ is $C$-special with respect to $p_1$, then
\[
\limsup_{k\to\infty}k_{\Hr^N}(Z_k,\rho(Z_K))\leq
\limsup_{k\to\infty}k_{\Hr^N}(Z_k,p_1(Z_K))+\limsup_{k\to\infty}k_{\Hr^N}(p_1(Z_k),\rho(Z_K))\leq
C,
\]
as stated.

On the other hand, if $\{Z_k\}$ is $p_1$-restricted (namely,
$\Re z_k\geq c \Im z_k$ for some $c>0$), from \eqref{rho} and
since $\Re z_k
>\|w_k\|^2$ it follows that $\{Z_k\}$ is also $\rho$-restricted
\end{proof}

\begin{remark}
It is worth to note explicitly that by Lemma
\ref{specialspecial} the condition  of being $C$-special and
that of being restricted do not depend on the chosen linear
projection.
\end{remark}

\begin{definition}
Let $h:\Hr^N\to \C$ be holomorphic. We say that $h$ has $E$-limit
$A\in \C$ at $\infty$, and we write
\[
{E\hbox{-}\lim}_{\Hr^N\ni (z,w)\to\infty}h(z,w)=A,
\]
if for any  sequence $\{Z_k\}\subset \Hr^N$ converging to $\infty$
which is $C$-special  for some $C\geq 0$ ($C$ depending on
$\{Z_k\}$) and  restricted, it follows that $\lim_{k\to
\infty}h(Z_k)=A$.

If the limit holds only for $0$-special, restricted sequences we
write
\[
{E^0\hbox{-}\lim}_{\Hr^N\ni (z,w)\to\infty}h(z,w)=A.
\]
\end{definition}

Next we state a version of the
Julia-Wolff-Carath\'eodory theorem due to Rudin for the unit
ball $\B^n$ (\cite[Thm. 8.5.6]{Ru}), using our
previous notations:

\begin{theorem}\label{JWC}
Let $\phi=(\phi_1,\phi'):\Hr^N\to\Hr^N$ be holomorphic, with
Denjoy-Wolff point $\infty$ and multiplier $\lambda\geq 1$. Let
$\rho:\Hr^N\to\Hr^N$ be a linear projection at $\infty$ and let
$\tr:\Hr^N\to \Hr$ be an associated left inverse. Then
\begin{enumerate}
    \item $\displaystyle{{E^0\hbox{-}\lim}_{\Hr^N\ni
    (z,w)\to\infty}\frac{\tr\circ\phi(z,w)}{\tr(z,w)}=\lambda}$,
    \item $\displaystyle{{E^0\hbox{-}\lim}_{\Hr^N\ni
    (z,w)\to\infty}\frac{\|\phi(z,w)-\rho\circ \phi(z,w)\|}{|\tr(z,w)|}=0}$.
\end{enumerate}
\end{theorem}

As a corollary we have the following:

\begin{lemma}\label{tech1}
Let $\phi=(\phi_1,\phi'):\Hr^N\to\Hr^N$ be holomorphic, with
Denjoy-Wolff point $\infty$ and multiplier $\lambda\geq 1$.
Assume  ${E\hbox{-}\lim}_{\Hr^N\ni
(z,w)\to\infty}\frac{\phi_1(z,w)}{z}$ exists. Then
\begin{enumerate}
    \item ${E\hbox{-}\lim}_{\Hr^N\ni
    (z,w)\to\infty}\frac{\phi_1(z,w)}{z}=\lambda$,
    \item ${E\hbox{-}\lim}_{\Hr^N\ni
    (z,w)\to\infty}\frac{\|\phi'(z,w)\|}{|z|}=0$.
\end{enumerate}
\end{lemma}
\begin{proof}
(1) It follows directly from Theorem \ref{JWC}.(1).

(2) Since $\phi(\Hr^N)\subseteq\Hr^N$ then $\Re \phi_1(z,w)\geq
\|\phi'(z,w)\|^2$ for all $(z,w)\in\Hr^N$. Thus dividing by $|z|^2$
and taking   limits, (2) follows from (1).
\end{proof}

The following technical proposition is needed in the proof of the Main Theorem.

\begin{proposition}\label{tech2}
Let $\phi=(\phi_1,\phi'):\Hr^N\to\Hr^N$ be holomorphic, with
Denjoy-Wolff point $\infty$ and multiplier $\lambda\geq 1$. Let
$\rho_0$ be a  linear projection at $\infty$ with a left
inverse $\tr_0$. Suppose the $E$-limit
${E\hbox{-}\lim}_{\Hr^N\ni
(z,w)\to\infty}\frac{\tr_0(v(z,w))}{\tr_0(z,w)}$ exists. Then
for any complex geodesic $f:\Hr\to \Hr^N$ with
$f(\infty)=\infty$ with left inverse $\tr_f:\Hr^N\to\Hr$ it
follows
\[{E\hbox{-}\lim}_{\Hr^N\ni (z,w)\to\infty}\frac{\tr_f\circ
\phi(z,w)}{\tr_f(z,w)}=\lambda.\]
\end{proposition}
\begin{proof}
Up to conjugation we can assume that $\rho_0=p_1$ and then by
hypothesis we know that the $E$-limit ${E\hbox{-}\lim}_{\Hr^N\ni
(z,w)\to\infty}\frac{\phi_1(z,w)}{z}$ exists and equals $\lambda$ by
Lemma \ref{tech1}. Given the complex geodesic $f$, there exists
$a\in \C^{N-1}$ and an automorphism $T\in\Aut(\Hr^N)$ of the type
\[
T(z,w)=(z+\|a\|^2+2\la w,a\ra, w+a)
\]
such that $T\circ f(\zeta)=(\zeta,0)$ and $\tr(z,w)=\pi_1 \circ
T(z,w)=z+\|a\|^2+2\la w, a\ra$, where, as usual $\pi_1(z,w)=z$. Thus
\[
\frac{\tr_f\circ \phi(z,w)}{\tr_f(z,w)}=\frac{\phi_1(z,w)+\|a\|^2+2\la
\phi'(z,w), a\ra}{z+\|a\|^2+2\la w,a\ra}=\frac{\phi_1(z,w)+\|a\|^2+2\la
\phi'(z,w), a\ra}{z(1+\|a\|^2/z+2\la w,a\ra/z)}.
\]
Taking into account that $1+\|a\|^2/z+2\la w,a\ra/z=1+o(|z|^{-1})$
since $|z|\geq \Re z\geq \|w\|^2$, the result follows from Lemma
\ref{tech1}.
\end{proof}

\end{document}